\documentclass[fleqn]{mat01}
\usepackage{times,mathtimy,amssymb,latexsym}
\begin{document}

\font\xx=msam5 at 10pt
\def\ab{\mbox{\xx{\char'03}}}

\font\bi=mtmib at 13.5pt
\def\f{\mbox{\bi{f}}}
\def\x{\mbox{\bi{x}}}
\def\u{\mbox{\bi{u}}}

\font\bc=mtmib at 10.4pt
\def\n{\mbox{\bc{n}}}

\font\bca=tib at 13.5pt
\def\?{\mbox{\bca{?}}}


\setcounter{page}{77}
\firstpage{77}

\def\defi{\trivlist\item[\hskip\labelsep{DEFINITION.}]}
\def\remark{\trivlist\item[\hskip\labelsep{\it Remark.}]}
\def\remarks{\trivlist\item[\hskip\labelsep{\it Remarks}]}
\def\noot{\trivlist\item[\hskip\labelsep{{\it Note.}}]}
\newtheorem{theo}{Theorem}[section]
\renewcommand\thetheo{\arabic{section}.\arabic{theo}}
\newtheorem{theor}[theo]{\bf Theorem}
\newtheorem{lem}{Lemma}

\title{When is  $\f( \x_1,\x_2, \ldots, \x_n )= \u_1(\x_1) + \u_2(\x_2)+  \cdots
+\u_n(\x_n) \?$}

\markboth{A K\l opotowski, M G Nadkarni and K P S Bhaskara Rao}{When is  $ f( x_1,x_2, \ldots, x_n )= u_1(x_1)+u_2(x_2)+  \cdots
+u_n(x_n) ?$}

\author{A K\L OPOTOWSKI$^{1}$, M G NADKARNI$^{2}$\\\noindent and K P S BHASKARA RAO$^{3}$}

\address{$^{1}$ Universit\'e Paris XIII,
 Institut Galil\'ee, 93430 Villetaneuse Cedex, France\\
\noindent $^{2}$ Department of Mathematics,
 University of Mumbai, Ka\-li\-na, Mumbai 400 098, India  \\
\noindent $^{3}$ Stat-Math Unit, Indian
 Statistical Institute, R.V.~College Post, Banglore~560 059, India\\
\noindent E-mail: klopot@math.univ-paris13.fr; nadkarni@math.mu.ac.in; kpsbrao@hotmail.com
}

\volume{113}

\mon{February}

\parts{1}

\Date{MS received 23 March 2001; revised 19 January 2002}

\begin{abstract}
We discuss  subsets  $S$ of ${\mathbb{R}}^n$ such 
that every real valued function $f$ on $S$  is of the form 
\begin{equation*}
f(x_1,x_2,\ldots,x_n) = u_1(x_1) + u_2(x_2) + \cdots + u_n(x_n),
\end{equation*}
and the related concepts and situations in analysis.
\end{abstract}


\maketitle

\section*{Introduction}

\noindent  Let $X_1,X_2,\ldots, X_n$ be non-empty sets. Let $S \subset
X_1\times X_2\times \cdots
\times X_n$. A point ${\underline{x}}\in S$ will look like ${\underline{x}}
=(x_1,x_2,\ldots, x_n)$.
Let $f: S \longrightarrow {\mathbb{R}}$ be a function. We say that $S$ is
{\it good  for} $f$, if we
can write $f$ in the form
\begin{equation*}
f({\underline{x}}) = u_1(x_1) + u_2(x_2) + \cdots + u_n(x_n),\quad
{\underline{x}}\in S,
\end{equation*}
where for each $i$, $u_i$ is a function from $X_i$ to ${\mathbb{R}}$. If
this holds for every function in a class \hbox{$\cal{A}$}\ \ of
functions on $S$, then we say that $S$ is {\it good for}
\hbox{$\cal{A}$}. We call $S$ {\it good}, if it is good for every
$f: S \longrightarrow \mathbb{R}$.


The purpose of this note is to give some descriptions of good sets and
comment on the connection of such sets with Kolmogorov's theorem on
superposition of functions and related questions in function algebras.
Connection with simplicial measures is also discussed (see \S5). For $n
= 2$ a geometric description of good sets is known, but this description
does not immediately generalize for the case $n >2$ (see \S4).

\section{Description of good sets}

Call a finite set $L =
\{{\underline{x}}^1,{\underline{x}}^2,\ldots,{\underline{x}}^k\}$ of
distinct points in $X_1\times X_2 \times \cdots \times X_n$ a {\it loop}
if:
\begin{enumerate}
\leftskip .5pc
\renewcommand{\labelenumi}{(\roman{enumi})}
\item there exist non-zero integers $p_1,p_2,\ldots,p_k$ such that
\begin{equation}
\hskip -.5cm p_1 {\underline {x}}^{1}+ p_2 {\underline {x}}^{2}+ \cdots + p_k
{\underline {x}}^{k} = 0,
\end{equation}
by which we mean that if $x_i^{j}$ is the $i{\rm{th}}$ coordinate of $
{\underline {x}}^{j}$, then
for each $i$, $1\leq i\leq n$, the formal sum $p_1x_i^{1} + p_2x_i^{2} +
\cdots + p_kx_i^{k} $
vanishes,

\item no proper subset of $L$ satisfies $(1)$.\\

Note that $(1)$ means that  $\sum_{j=1}^kp_j{\bf 1}_{\{x^j_i\}} = 0$  for
each $i$.
\end{enumerate}

\begin{remark}
For $n=2$ the integers $p_i$ can be chosen to be $+1$ or $-1$, but for
$n \geq 3$ this fails and there is no universal upper bound (depending
on $n$) on the integers $p_1, p_2,\ldots, p_k$ (see \S4).
\end{remark}

\begin{theor}[\!] 
Let $S \subset X_1 \times X_2 \times \cdots \times X_n$ and let $f : S
\longrightarrow {\mathbb{R}}$ be such that whenever the formal sum
$\sum_{j=1}^k p_j \underline{x}^{(j)}=0${\rm ,} then $\sum_{j=1}^k p_j
f(\underline{x}^{(j)})=0.$ Then there exist real valued functions $u_1,
u_2, \ldots, u_n$ defined on $X_1, X_2,\ldots , X_n$ respectively such
that 
\begin{equation}
 f(\underline{x}) = f (x_1, x_2, \ldots, x_n) = u_1(x_1)+u_2(x_2)
+ \cdots+ u_n(x_n),
\end{equation}
for all $(x_1,x_2,\ldots,x_n) \in S$.
\end{theor}

\begin{proof}
It is clear that if $f$ is of the form $(2)$, then for any loop $L
=\{{\underline {x}}^{1}, {\underline {x}^{2}},\ldots,{\underline
{x}}^{k}\}$ of points in $S$ the sum $\sum_{j=1}^kp_jf(\underline{x}^j)$
vanishes.

Assume now that for any loop $L =\{{\underline {x}}^{1},
{\underline {x}^{2}},\ldots,{\underline {x}}^{k}\}$ of points in $S$ the
sum $\sum_{j=1}^kp_jf(\underline{x}^j)$ vanishes. We can suppose without
loss of generality that $X_i\cap X_j = \emptyset$ for $i\neq j$. Let
$\Omega = X_1 \cup X_2 \cup \cdots\cup X_n$. Every $\underline{x}=
(x_1,x_2,\ldots,x_n) \in S$ has associated to it a subset of $\Omega$,
namely the set $\{x_1,x_2,\ldots,x_n\}$ with $n$ points. Let 
\begin{equation*}
{\cal C} = \{\{x_1,x_2,\ldots,x_n\}:\ (x_1,x_2,\ldots,x_n)\in S\}.
\end{equation*}
Then ${\cal C}$ is a collection of subsets of $\Omega$. Define on
$\cal{C}$ the function $\mu$ by 
\begin{equation*}
\mu(\{x_1,x_2,\ldots,x_n\})= f(x_1,x_2,\ldots,x_n).
\end{equation*}
 The class $\cal{V}$ of functions of the form
$\sum_{j=1}^lr_j{\bf 1}_{C_j}$, $r_j$ rational, $C_j \in \cal{C}$, $l
\geq 1$, is a vector space over the field of rational numbers and the
condition that for any loop $L =\{{\underline {x}}^{1}, {\underline
{x}^{2}},\ldots,{\underline {x}}^{k}\}$ of points in $S$ the sum
$\sum_{j=1}^kp_jf(\underline{x}^j)$ vanishes, ensures that the map $T$
on $\cal{V}$ defined by 
\begin{equation*}
T\left(\sum_{j=1}^lr_j{\bf 1}_{C_j}\right) =\sum_{j=1}^lr_j\mu(C_j)
\end{equation*}
is well defined and linear. We extend this map linearly to the larger
class $\cal{W}$ of functions of the form $\sum_{j=1}^lr_j{\bf 1}_{C_j}$,
$r_j$ rational, $C_j \subset \Omega,\ l \geq 1$, and continue to denote
the extended map by $T$. Let us define $u_i : X_i \longrightarrow
{\mathbb{R}}$ by $u_i(x_i)= T{\bf 1}_{\{ x_i\}}$ for $x_i\in X_i$, $1
\leq i\leq n$. Now, for any $\underline{x}= (x_1, x_2, \ldots x_n) \in
S$, 
\begin{align*}
f(\underline{x}) &= \mu (\{x_1, x_2, \ldots , x_n\}) = T{\bf1}_{\{
x_1, x_2, \ldots, x_n\}}\\
&= T{\bf 1}_{\{ x_1\}} +\cdots + T{\bf
1}_{\{x_n\}}= u_1(x_1) + u_2(x_2) + \cdots + u_n(x_n).\\[-1.3cm]
\end{align*}\vspace{.5cm}
\hfill $\ab$
\end{proof}


\begin{theor}[\!]
A set $S \subset X_1\times X_2 \times\cdots \times X_n$ is good if and
only if $S$ has no loop in it.
\end{theor}
\pagebreak

\begin{proof}
If $S \subset X_1\times X_2\times\cdots\times X_n$ does not admit a
loop, then the hypothesis of Theorem~1.1 is vacuously satisfied
and so any real valued function on $S$ is of the form
\begin{equation*}
f(x_1,x_2,\ldots,x_n) = u_1(x_1) + u_2(x_2)+\cdots + u_n(x_n),
\end{equation*}
$(x_1,x_2,\ldots,x_n) \in S$, where $u_1, u_2,\ldots, u_n$ are functions
defined on $X_1, X_2,\ldots, X_n$ respectively. On the other hand, if
$S$ admits a loop then an $f$ violating the condition of Theorem 1.1 can
be constructed easily, so Theorem~1.2 follows.\hfill $\ab$
\end{proof}

\begin{remarks}
$\left.\right.$\vspace{.2pc}

\begin{enumerate}
\leftskip .5pc
\renewcommand{\labelenumi}{(\roman{enumi})}
\item Clearly Theorem~1.2 is also valid for complex-valued functions
$f$. One simply treats real and imaginary parts separately. In the
sequel we shall take $f$ to be complex valued.

\item If $S \subset {\mathbb{R}^n}$ is good and the canonical projections of
$S$ on the coordinate axes are pairwise disjoint, then clearly we can
choose the $u_i$'s all equal. If $S \subset {\mathbb{R}^n}$ is good,
then for any $\underline{c} \in \mathbb{R}^n$ the set $S +\underline{c}$
is also good and, when $S$ is bounded, for a suitable $\underline{c}$ the canonical
projections of $S + \underline{c}$ on the coordinate axes are pairwise
disjoint, so one can choose the functions $u_i$, for a given $f$ on such
an $S + {\underline {c}}$, to be the same.\vspace{.5pc}
\end{enumerate}
\end{remarks}
To end this section we shall give a description of good subsets $S$ of
$X_1\times X_2 \times \cdots \times X_n$, when all the sets
$X_1,X_2,\ldots , X_n$ are finite, i.e, $\hbox{card}\, X_i=
m_i<+\infty,\ 1\le i\le n.$
 
Let $\Pi_i:X_1 \times X_2 \times \cdots \times X_n\longrightarrow X_i$,
$1\le i\le n$, be the canonical projections on $X_i$. If $S$ is good,
then any function $f:S\longrightarrow {\mathbb R},\ f =u_1 +u_2 +\cdots
+u_n,$ is completely deter\-mi\-ned by the values of $u_i$ on $\Pi_i S$,
$1\le i\le n$. Hence we can assume in addition that $\Pi_i S=X_i$, $1\le
i\le n$.

Let $X_i=\{x_1^{(i)},x_2^{(i)},\ldots,x_{m_i}^{(i)}\},\ 1\le i\le n,$
and $S=\{s_1,s_2,\ldots,s_k\},$ where
\begin{equation*}
s_j=(x_{j_1}^{(1)},x_{j_2}^{(2)},\ldots ,x_{j_n}^{(n)})\qquad 1\le j\le k,\qquad
1 \leq j_i \leq m_i.
\end{equation*}

We consider the $k\times (m_1+m_2+\cdots +m_n)$-matrix $M$ (called {\it
the matrix of} $S$) with rows $M_j,\ 1\le j\le k,$ given by
\begin{equation*}
M_j=(0,\ldots,0,1,0,\ldots,0,1,0,\ldots ,0,1,0,\ldots , 0),
\end{equation*}
where $1$ occurs at the places $j_1$, $m_1+j_2 $, $m_1+m_2 +j_3 $, etc.
corresponding to the subscripts in the point
$s_j=(x_{j_1}^{(1)},x_{j_2}^{(2)},\ldots,x_{j_n}^{(n)}),\ 1\le j\le k.$
Since $S$ is good, 
\begin{align*}
f(s_j) &=f(x_{j_1}^{(1)},x_{j_2}^{(2)},\ldots,x_{j_n}^{(n)}) \\
&= u_1(x_{j_1}^{(1)}) + u_2(x_{j_2}^{(2)})+\cdots +
u_n(x_{j_n}^{(n)}),\quad 1\leq j \leq k.
\end{align*}
We put 
\begin{align*}
u_1(x_{1}^{(1)})&=\alpha_{1}^{(1)},
\ldots,u_1(x_{m_1}^{(1)})=\alpha_{m_1}^{(1)},\\
u_2(x_{1}^{(2)})&=\alpha_{ 1}^{(2)},
\ldots,u_2(x_{m_2}^{(2)})=\alpha_{m_2}^{(2)},\\
&\qquad\qquad\cdots\\
u_n(x_{1}^{(n)})&=\alpha_{ 1}^{(n)},
\ldots,u_n(x_{m_n}^{(n)})=\alpha_{m_n}^{(n)}.
\end{align*}
The relation $(2)$ gives us $k$ equalities
\begin{equation*}
\alpha_{j_1}^{(1)} +\alpha_{j_2}^{(2)}+\cdots
+\alpha_{j_n}^{(n)}=f( s_j),\quad 1\le j\le k.
\end{equation*} 
In other words, the column vector
\begin{equation*}
(\alpha_{1}^{(1)}, \ldots , \alpha_{m_1}^{(1)},\alpha_{
1}^{(2)},\ldots, \alpha_{m_2}^{(2)},\ldots ,\alpha_{ 1}^{(n)}, \ldots ,
\alpha_{m_n}^{(n)})^t\in {\mathbb R}^{m_1+m_2+\cdots+m_n}
\end{equation*}
is a solution of the matrix equation
\begin{equation}
M\vec \alpha =\vec z,
\end{equation}
where $\vec z= (f( s_1), f( s_2), \ldots, f( s_k))^t\in {\mathbb
R}^{k}.$

Since $S$ is good, we know that $(3)$ has solution
for every $\vec z $. Since $M$ has $m_1+m_2+\cdots+m_n$ columns and
since the $n-1$ vectors
\begin{equation*}
\begin{array}{l}
(\underbrace{1,1,\ldots,1}_{m_1\
{\rm times}},\underbrace{-1,\ldots,-1 }_{m_2\ {\rm times}},0,0,0,\ldots)^t\\[2pc]
(\underbrace{1,1,\ldots,1}_{m_1\ {\rm times}},\underbrace{0,\ldots,0}_{m_2\
{\rm times}},\underbrace{-1,\ldots,-1 }_{m_3\ {\rm times}},0,0,0,\ldots)^t\\[1.5pc]
\cdots\\[.5pc]
(\underbrace{1,1,\ldots,1}_{m_1\ {\rm times}},0,\ldots,
0,\underbrace{-1,\ldots,-1 }_{m_n\ {\rm times}})^t
\end{array}
\end{equation*}
are linearly independent
solutions of the homogeneous equation $M\vec\alpha= \vec 0$, we see that
the rank of $M$ is at most $m_1+m_2+\cdots+m_n-(n-1)$. Clearly $k$ 
cannot exceed the rank of $M$. On the other hand the union of $n$ sets
\begin{align*}
&(X_1\times \{x_2\}\times \cdots \times \{x_n\}) \cup (\{x_1\}\times
X_2\times\{x_3\}\times \cdots\\
&\quad\,\,\, \times \{x_n\})\cup\cdots \cup
(\{x_1\}\times \cdots \times \{x_{n-1}\}\times X_n)
\end{align*}
is a good subset of $X_1 \times X_2 \times \cdots \times X_n$ of
cardinality $m_1+m_2 +\cdots+m_n-(n-1)$. It is clear that if the rank of
$M$ is $k$ and $k\leq m_{1} + m_{2} + \cdots + m_{n} - (n - 1)$ then $S$
is good. We have proved:

\begin{theor}[\!]
Let $S$ be a finite subset of $ X_1\times X_2 \times\cdots \times X_n$
of cardinality $k$ and let $m_i$ denote the cardinality of $\Pi_iS${\rm ,} the
canonical projection of $S$ on $X_i$. Then $S$ is good if and only if
$k\leq m_1+m_2+\cdots+m_n-(n-1) $ and the matrix $M$ of $S$ defined
above has \hbox{\rm rank} $k$. There always exist a good set of cardinality $k \leq
m_1+m_2+\cdots+m_n-(n-1) $.
\end{theor}
Let us remark also that the procedure described in  Proposition~2.7 of \cite{ckn} does
not work even in the three-dimensional case.\\

\section{Sequentially good sets}

We say that $S$ is {\it sequentially good for a complex valued function}
$f$ defined on $S$ if
\begin{equation*}
f(x_1, x_2,\ldots,x_n) = \lim_{k\to \infty}(u_{1,k}(x_1) + u_{2,k}(x_2)+
\cdots + u_{n,k}(x_n)),
\end{equation*}
where $(x_1,x_2,\ldots, x_n) \in S$ and $u_{1,k}, u_{2,k}, \ldots,
u_{n,k},\ k = 1,2,3,\ldots$ are functions on $X_1, X_2, \ldots, X_n$
respectively. If $S$ is sequentially good for every function on $S$,
then we say that $S$ is {\it sequentially good}. It is clear that if a
set $S$ is good for $f$, then it is sequentially good for $f$. The
converse holds in view of Theorem~1.2. Indeed, if $S$ is
sequentially good for $f$, but not good for $f$, then there exists a
loop $L =\{{\underline {x}}^{1}, {\underline {x}^{2}},\ldots,{\underline
{x}}^{k}\}$ of points in $S$ such that the sum
$\sum_{j=1}^kp_jf(\underline{x}^j)$ does not vanish, and at the same
time $f$ is the pointwise limit of a sequence of functions $g_n,\ n =
1,2,\ldots$ such that for each $g_n$,
$\sum_{j=1}^kp_jg_n(\underline{x}^j)$ vanishes. The contradiction shows
that $S$ is good for $f$.

Say that a subset $S$ of $X_1\times X_2 \times\cdots\times X_n$ is {\it
sequentially good for a collection} ${\cal{F}}$ of functions on $S$, if
every $f\in {\cal{F}}$ is of the form
\begin{equation*}
f(x_1,x_2,\ldots,x_n) = \lim_{k\to \infty}(u_{1,k}(x_1)+u_{2,k}(x_2)
+\cdots + u_{n,k}(x_n)),
\end{equation*}
$(x_1,x_2,\ldots, x_n) \in S ,\ u_{1,k}, u_{2,k},\ldots, u_{n,k},\ k =
1,2,\ldots $ being functions on $X_1,X_2,\ldots,$ $X_n$ respectively.

Assume now that $S$ is sequentially good for an algebra ${\cal{F}}$ of
functions on $S$ which is closed under conjugation, separates points and
contains constants. Then in fact $S$ is sequentially good (hence good).
For otherwise $S$ will admit a loop $L$. The restriction of functions in
${\cal{F}}\,\,\!$ to $L$ (denoted by ${\cal{F}}\!\!\!\mid_L$) is an algebra of
functions on $L$, closed under conjugation, separating points and
containing constants. Since $L$ is a finite set (hence compact in the
discrete topology), by Stone--Weierstrass theorem, the algebra
${\cal{F}}\!\!\mid_L$ is dense in the collection of functions on $L$, hence
actually equal to the collection of all functions on the finite set $L$.
Since $L$ is sequentially good for all functions on $L$, we see by our
earlier conclusion that $L$ is good and so not a loop. The contradiction
shows that $S$ is good. We have proved:

\begin{theor}[\!]
The following are equivalent for a set
 $S\subset X_1\times X_2 \times \cdots \times X_n${\rm :}
\begin{enumerate}
\leftskip .5pc
\renewcommand{\labelenumi}{\rm (\roman{enumi})}
\item  $S$ is good{\rm ,}

\item  $S$ is sequentially good{\rm ,}

\item every finite subset of $S$ is good{\rm ,}

\item $S$ is sequentially good for an algebra of functions
on $S${\rm ,} which is closed under conjugation{\rm ,} separates points of $S$ and
contains constants.
\end{enumerate}
\end{theor}

\section{Sequentially good measures}

Let $X_1,X_2, \ldots ,X_n$ be Polish spaces. Call a probability measure
$\mu$ on Borel subsets of $ \Omega = X_1 \times X_2 \times \cdots \times
X_n$ {\it sequentially good for a collection} ${\cal{F}}$ of
complex-valued functions on $ \Omega$ if every function $f \in
{\cal{F}}$ is of the form
\begin{align*}
f(x_1,x_2,\ldots, x_n) =
\lim_{k\to\infty}(u_{1,k}(x_1) + u_{2,k}(x_2) +\cdots
 +
u_{n,k}(x_n)),\ \ \mu - \ \hbox{a.e.},
\end{align*}
where $u_{1,k},u_{2,k},\ldots, u_{n,k},\ k = 1,2,\ldots$ are Borel
measurable.\vspace{.3pc}

Let $ A_1,A_2,A_3, \ldots$ be a countable collection of Borel subsets of
$\Omega$ which is closed under finite unions and compliments and
separates points of $\Omega$. Let $\mu$ be a sequentially good
probability measure for the countable collection of functions ${\bf
1}_{A_i},\ i = 1,2,3,\ldots$. Then there is a Borel subset $S$ of full
$\mu$ measure which is sequentially good for the collection ${\bf
1}_{A_i}, i = 1,2,3,\ldots$. The set $S$ continues to be sequentially
good for the algebra ${\cal{A}}$ of finite linear combinations of ${\bf
1}_{A_i},\ i = 1,2,3,\ldots $with complex coefficients, an algebra which
is closed under conjugation, separates points and contains constants. By
Theorem~2.1 the set $S$ is sequentially good, hence a good set.
We have proved:

\begin{theor}[\!]
If $\mu$ is sequentially good for the countable collection of indicator
functions ${\bf 1}_{A_i},\ i~=~1,2,3,\ldots$ of sets in a countable
field of Borel sets which separate points of $X_1\times X_2\times
\cdots\times X_n${\rm ,} then $\mu$ admits a Borel support $S$ which is
good.
\end{theor}

\section{Cases ${\n} ={\bf 2}$ and ${\n}>{\bf 2}$}

A good subset of ${\mathbb{R}}^2$ has a geometric description which does
not seem to be available for $n>2$.

Two arbitrary points $(x,y)$, $(z,w)$ in $S\subseteq X\times Y $ ($S$ is
not ne\-cessarily good) are said to be {\it linked} (and we write
$(x,y)L(z,w)$), if there exists a finite sequence of points $(x_1,y_1)$,
$(x_2,y_2)$, $\ldots$, $(x_n,y_n)$ in $S$ (called {\it a link} of {\it
length} $n$ joining $(x,y)$ to $(z,w) $) such that:

\begin{enumerate}
\leftskip .5pc
\renewcommand{\labelenumi}{(\roman{enumi})}
\item $(x_1,y_1) =(x,y),\ (x_n,y_n)=(z,w);$

\item for any $i, 1\le i\le n-1 $ {\it exactly one} of the
following equalities holds:
\begin{equation*}
\hskip -.5cm x_i=x_{i+1},\quad   y_i=y_{i+1};
\end{equation*}
\item for any $i$, $1\leq i\leq n-2$, it is not possible to have
$x_i=x_{i+1}=x_{i+2}$ or  $ y_i=y_{i+1}=y_{i+2}$.
\end{enumerate}\vspace{-.4pc}

Note that $L$ is an equivalence relation. An equivalence class of $L$ is
called {\it a linked component} of $S$. If $(x,y)\in S$, then the
equivalence class to which $(x,y)$ belongs is called the {\it linked
component} of $(x,y)$. Two points $(x,y),(z,w)\in S $ are said to be
{\it uniquely linked}, if there is a unique link joining $(x ,y )$ to
$(z,w)$. A linked component of $S\subseteq X\times Y $ is said to be
{\it uniquely linked} if any two points in it are uniquely linked.

One can prove (see \cite{ckn,HW}) that a subset $S \subset
X\times Y$ is good if and only if each of its linked components is
uniquely linked. See \cite{kn,KN} for more discussion on good sets for $n = 2$.

A geometric description of good subsets $S$ of $X \times Y \times Z$,
and more generally of $X_1 \times X_2 \times \cdots \times X_n$ is not
available. We only have a partial answer. We consider here the case $n =
3$. For $n>3$ the notion of a link and linked component can be similarly
defined.

\begin{defi}$\left.\right.$\vspace{.3pc}

\noindent Two arbitrary points $(x,y,z),(p,q,r)\in S\subseteq X\times Y\times Z$
are said to be {\it linked} (and we write $(x,y,z)L(p,q,r)$), if there
exists a finite sequence of points
$\{(x_1,y_1,z_1),(x_2,y_2,z_2),\ldots,$ $(x_n,y_n,z_n)\}$ in $S$ (called a
{\it link} joining $(x,y,z)$ to $(p,q,r) $) such that:

\begin{enumerate}
\leftskip .5pc
\renewcommand{\labelenumi}{(\roman{enumi})}
\item  $(x_1,y_1,z_1) =(x,y,z),\ (x_n,y_n,z_n)=(p,q,r)$,

\item  for any $1\le i\le n-1 $ {\it exactly one} of the
following holds
\begin{equation*}
\hskip -.5cm x_i\neq x_{i+1},\quad   y_i\neq y_{i+1},\quad   z_i\neq z_{i+1},
\end{equation*}
\item for any $i$, $1\leq i\leq n-2$, none of the following
holds:
\begin{align*}
\hskip -.5cm &(x_i\neq x_{i+1} \quad \hbox{and}\quad  x_{i+1}\neq x_{i+2}),\\
\hskip -.5cm &(y_i\neq y_{i+1} \quad \hbox{and}\quad y_{i+1}\neq y_{i+2}),\\
\hskip -.5cm &(z_i\neq z_{i+1} \quad \hbox{and}\quad z_{i+1}\neq z_{i+2}).
\end{align*}
\end{enumerate}
\end{defi}

As before $L$ is an equivalence relation. A uniquely linked set is
similarly defined. An equivalence class of $L$ is called {\it a linked
component} of $S$. We call $S$ {\it linked}, if it has only one linked
component. As in the case of two-dimensional sets, one can prove:\vspace{.3pc}

{\it A linked set $S \subset X \times  Y\times Z$ is good if and only if it
is uniquely linked.}\vspace{.3pc}

However, it is not true that a subset $S\subset {\mathbb{R}}^3$ is good if
each linked component is
uniquely linked, as the following example shows:

The set $ \{ (0,0,0),(0,0,1),(1,1,0),(1,1,1)\}$ has two uniquely linked
components, namely, $\{(0,0,0), (0,0,1)\}$ and $\{ (1,1,0), (1,1,1)\}$,
but it is not a good set, as can be seen by writing four linear
equations in six unknowns $u(0), v(0), w(0), u(1), v(1), w(1)$ 
onto~${\mathbb{R}}^4$.

In case $n = 2$, the coefficients $p_i$ in the definition of a loop can
be chosen to be $+1$ or $-1$. However, for $n>2$ the coefficients $p_i$
do not have a universal bound (depending only on $n$). Here are two
examples: The set
\begin{equation*}
\{(0,0,0),(0,0,1),(0,1,0),(1,0,0),(1,1,1)\}
\end{equation*}
is not a good subset of $\{0,1\}^3$. It is also a loop, because the
formal sum
\begin{equation*}
2(0,0,0)-(0,0,1)-(0,1,0)-(1,0,0)+(1,1,1)
\end{equation*}
is equal to $0$. This loop is minimal (i.e. each of its proper subset is
good) and one cannot have the formal sum above vanish with all the
coefficients equal to $+1$ or $ -1$. For the second example, let $X_1 =
X_2 = X_3 = {\mathbb{R}}$. For the obvious loop described by the
following expression not all $p_i's$ can be chosen less than five.

\begin{equation*}
\begin{array}{l@{\qquad}l@{\qquad}l}
 5~(1~1~1)\\
 -(2~3~ 1)&-(12, 1, 13)&-(1, 22, 23)\\
 -(4~ 5~ 1) &-(14, 1, 15)& -(1, 24, 25) \\
 -(6, 7, 1) &-(16, 1, 17) &-(1, 26, 27) \\
 -(8, 9, 1) &-(18, 1, 19) &-(1, 28, 29)\\
 -(10, 11, 1) &-(20, 1, 21) &-(1, 30, 31)\\
 +(2, 5, 13) &+(12, 22, 25)&\\
 +(4, 7, 15) &+(14, 24, 27)&\\
 +(6, 9, 17) &+(16, 26, 29)&\\
 +(8, 11, 19) &+(18, 28, 31)&\\
 +(10, 3, 21) &+(20, 30, 23)&
 \end{array}
\end{equation*}

The above example can be modified so that at least one $p_i$ is bigger
than $P$, a pre-assigned positive integer $ \geq2$.

\section{Discussions}

As a solution to Hilbert's 13th problem, Kolmogorov (see \cite{Kh,Ko,Lo}) proved that one can imbed the unit cube $E^n =
[0,1]^n$ in $\mathbb{R}^{2n +1}$ homeomorphically by a map of the type
$\psi: (x_1,\ldots, x_n) \longrightarrow
(\sum_{p=1}^n\psi_{1,p}(x_p),\ldots, \sum_{p=1}^{n}\psi_{2n+1,p}(x_p))$,
with $\psi_{q,p}$ continuous and monotonic increasing on $[0,1]$, such
that every continuous function $g$ on $\psi(E^n)$ is of the form
\begin{equation*}
g(y_1,\ldots,y_{2n+1}) =\sum_{q=1}^{2n+1}g_q(y_q).
\end{equation*}
In particular this implies that $\psi(E^n)$ is a good set for complex
valued continuous functions, and since such functions form an algebra
closed under conjugation, contain constants, and separate points, we see
by Theorem~2.1 that $\psi (E^n)$ is a good set. It has been
observed by Lorentz~\cite{Lo} that $\psi$ can be chosen so that
$g_1,\ldots, g_n$ are all equal. Remark (ii) following 
Theorem~1.2 shows how this may be arranged.

Two questions naturally arise:

\noindent (A) describe compact subsets of $C \subset \mathbb{R}^n$
such that every continuous function $g$ on $C$ is of the form 
\begin{equation*}
g(y_1,\ldots,y_{n}) =\sum_{q=1}^{n}g_q(y_q),
\end{equation*}
with $g_1,\ldots,g_n$ continuous,

\noindent (B) describe compact subsets of $C \subset \mathbb{R}^n$
such that every continuous function $g$ on $C$ is of the form
\begin{equation*}
g(y_1,\ldots,y_{n}) = \lim_{l\to \infty}\sum_{q=1}^{n}g_{q,l}(y_q),
\end{equation*}
with $g_{q,l},\ 1 \leq q \leq n,\  l = 1,2,\ldots$ continuous.\vspace{.3pc}

For $n = 2$ these questions are well discussed in the literature. For
question (A) a necessary and sufficient condition on $C$ is that it be
loopfree (i.e., a good set) and the lengths of links in $C$ be bounded
\cite{Ma,rdm,rdm1}. For question (B) a sufficient
condition is that $C$ be loopfree and that linked components be closed
\cite{Ma1} or more generally that linked components admit a Borel
cross-section \cite{KN1}.

For $n>2$ natural analogues of these are not known since a good
definition of linked component is not available (see also \cite{spr,str}). Theorem~2.1 however shows that a necessary condition
on $C$ for both question (A) and (B) is that $C$ be loopfree.\vspace{.3pc}

Let $X_1, X_2, \ldots, X_n$ be Polish spaces and let $\Omega = X_1\times
X_2\times\cdots \times X_n$. A probability measure $\mu$ on $\Omega $ is
said to be {\it simplicial}, if $\mu$ is an extreme point of the convex
set of all probability measures $\lambda$ on $\Omega$, whose 
one-dimensional marginals are the same as those of $\mu$. Let $\mu$ be a
simplicial measure and let $\mu_1, \mu_2,\ldots, \mu_n$ denote the
one-dimensional marginals of $\mu$. A theorem of Lindenstrauss
\cite{Li} and Douglas \cite{D} states that:\vspace{.5pc}

{\it A probability measure $\mu$ on $\Omega$ is simplicial if and only
if the collection of functions of the form 
\begin{equation*}
f(x_1,x_2,\ldots, x_n) = u_1(x_1)+u_2(x_2)+\cdots + u_n(x_n),
\end{equation*}
where $u_i\in L^1(X_i,\mu_i),\ 1\leq i\leq n${\rm ,} is dense in
$L^1(\Omega,\mu)$.}\vspace{.5pc}

This theorem is usually proved for $n = 2$, but the same proof holds for
any $n$. It is clear from this theorem that a simplicial measure is 
sequentially good for the functions ${\bf 1}_{A_i},\ i = 1,2,3,\ldots$,
where $\{A_i: i = 1,2,3,\ldots\}$ form a countable field of Borel sets
which separate points of $\Omega$ and so by Theorem~3.1 admits a
Borel support which is a good set. We have proved:

\begin{theor}[\!]
A simplicial measure admits a good Borel set as support.
\end{theor}

For $n=2$ this result is due to Bene\v{s} and
\v{S}t\v{e}p\'{a}n \cite{BS,bs}.

If $\lambda_1, \lambda_2,\cdots, \lambda_n$ are continuous probability
measures on $X_1,X_2,\ldots,X_n$ respectively, then it is an easy
consequence of Fubini theorem that any Borel set of positive
$\lambda_1\times\lambda_2\times\cdots\times \lambda_n$ measure contains
a loop of the type $B_1\times B_2\times \cdots \times B_n$ with each
$B_i$ a two point set. Since a simplicial measure admits a good Borel
set as support, we see that {\it a simplicial measure is singular to
$\lambda_1\times\lambda_2\times\cdots\times \lambda_n$ for any choice of
continuous $\lambda_1, \lambda_2,\ldots, \lambda_n$ on
$X_1,X_2,\ldots,X_n$ respectively} (see \cite{Li,ST} for the
case $n = 2$).

Let us briefly return to question (B) above and let $C$ be a compact
subset of $\mathbb{R}^n$ such that every continuous function on $C$ is
approximable as described there. Then every probability measure on $C$
is simplicial. For, if $\mu_1$ and $\mu_2$ are two distinct probability
measures on Borel subsets of $C$ with the same one-dimensional marginals
then $\mu_1 - \mu_2$ is a non-trivial signed measure which integrates
all continuous functions on $C$ to zero, which is not possible.

\begin{remark}
For a discussion of Hilbert's 13th problem from 
algebraic point of view see \cite{A,AS}.  
\end{remark}

\end{document}